\documentclass[preprint,review,10pt]{elsarticle}

\usepackage{amssymb}
\RequirePackage{graphicx}
\journal{}

\begin{document}

\begin{frontmatter}


\author{Ping Sun}
%\fntext[label2]{Tel.:+86 024-83686202}
\ead{plsun@mail.neu.edu.cn}

\title{A probabilistic approach for enumeration of certain Young tableaux}


%\author{}

\address{Department of Mathematics, Northeastern University, Shenyang, 110004, China}

\begin{abstract}

In this paper we establish an order statistics model of Young tableaux. Multiple integration over nested simplexes is applied to the enumeration of Young tableaux. A brief proof of Frobenius-Young's and Aitken's formulas is given. Partially standard Young tableaux and special truncated shapes including tableaux with a hole are discussed, the associated product formulas are given.

\end{abstract}

\begin{keyword}

Young tableaux \sep hook-length formulas \sep order statistics \sep determinants \sep multiple integrals

\MSC[2010] 05D40 \sep 05E10 \sep 05A17 \sep 62G30 \sep 26B15

\end{keyword}

\end{frontmatter}

%\linenumbers


\section{Introduction}

Probabilistic method involving techniques of continuous mathematics in combinatorial enumeration problem is always interesting. This paper shows that order statistics may play an important role. Suppose that $\xi_{1,n} \leq \xi_{2,n}\leq \cdots \leq \xi_{n,n}$ are order statistics from independent and identically distributed (i.i.d) random variables with uniform distribution on $[0,1]$. In a previous paper \cite{SUW}, the uniform distribution on $[0,1]$ was applied to Stirling numbers of the second kind $S(n,k)$, we have obtained that
$$
S(n,k) = {n \choose k} {\bf E} (\xi_{1,k} + \cdots + \xi_{k,k})^{n-k}, \eqno(1.1)
$$
here {\bf E} is mathematical expectation operator.

There is a generalization of (1.1) of the following form with $m \geq k$:
$$
S(n,k) = {n+m-k \choose m} {\bf E} (\xi_{1,m} + \cdots + \xi_{k,m})^{n-k}. \eqno(1.2)
$$
Clearly (1.2) is equivalent to (see (2.3) below)
$$
S(n,k) = (n+m-k)_m \mathop {\int {\int {...\int {} } } }\limits_{\tiny
   {0 \leq x_1 < \cdots < x_m \leq 1} } (x_1+x_2+\cdots+x_k)^{n-k} dx_1  \cdots dx_m
$$
$$
=(n+m-k)_m \mathop {\int {\int {...\int {} } } }\limits_{\tiny \begin{array}{*{10}c}
   x_i \geq 0, 1 \leq i \leq m  \\
   {x_1 + \cdots + x_m \leq 1}  \\
\end{array} } (x_1+2x_2+\cdots+kx_k)^{n-k} dx_1  \cdots dx_m. \eqno(1.3)
$$
Here $(n)_k=n(n-1)\cdots(n-k+1)$. In fact, one may verify (1.3) by using of Dirichlet's multiple integral.

Recently, the author applied order statistics to multiple integrals, classical results of Dirichlet's and Liouville's integrals were given new proofs \cite{SU}. In this paper we develop further applications. Certain nested groups of order statistics are exposed to corresponding with Young tableaux, so it is possible to obtain the numbers of Young tableaux by the method of multiple integration.

This paper is organized as follows. In Section 2 we review the relevant definitions and notation. In Section 3 we introduce the probabilistic model and explain how should it be applied to the enumeration of Young tableaux. In Sections 4 we state the main result and find a constructive proof of two classical formulas for (skew) standard Young tableaux. Section 5 gives the product formulas for the number of tableaux of specific shapes: partially standard Young tableaux, three variations of Young tableaux with three rows removed by one box; illustrated as

\vspace{5pt}

\begin{picture}(50,50)(-20,0)
\put(0,0){\line(10,0){60}} \put(0,10){\line(10,0){60}} \put(0,20){\line(10,0){50}} \put(0,30){\line(10,0){40}} \put(0,40){\line(10,0){30}}
\put(0,0){\line(0,10){39.5}} \put(10,0){\line(0,10){39.5}} \put(20,0){\line(0,10){39.5}} \put(30,0){\line(0,10){39.5}}
\put(40,0){\line(0,10){29.5}} \put(50,0){\line(0,10){19.5}} \put(60,0){\line(0,10){9.5}} \put(-30,-20){{\footnotesize $f_{(m+1,\cdots,m+k) \subset (m+r+1,\cdots,m+r+k)}$}}

\put(110,10){\line(10,0){20}} \put(110,20){\line(10,0){50}} \put(110,30){\line(10,0){50}} \put(110,40){\line(10,0){40}}
\put(110,10){\line(0,10){29.5}} \put(120,10){\line(0,10){29.5}} \put(130,10){\line(0,10){29.5}} \put(140,20){\line(0,10){19.5}}
\put(150,20){\line(0,10){19.5}} \put(160,20){\line(0,10){10}} \put(110,-10){{\footnotesize $f_{m,m+1,m-k}$}}

\put(180,10){\line(10,0){40}} \put(180,20){\line(10,0){40}} \put(180,30){\line(10,0){60}} \put(180,40){\line(10,0){60}}
\put(180,10){\line(0,10){29.5}} \put(190,10){\line(0,10){29.5}} \put(200,10){\line(0,10){29.5}} \put(210,10){\line(0,10){29.5}}
\put(220.5,10){\line(0,10){10}} \put(220,30){\line(0,10){9.5}} \put(230,30){\line(0,10){9.5}} \put(240,30){\line(0,10){9.5}} \put(180,-10){{\footnotesize $f_{m+k,m-1,m}$}}

\put(260,10){\line(10,0){30}} \put(260,20){\line(10,0){29.5}} \put(260,30){\line(10,0){70}} \put(260,40){\line(10,0){70}}
\put(260,10){\line(0,10){29.5}} \put(270,10){\line(0,10){29.5}} \put(280,10){\line(0,10){29.5}} \put(290,10){\line(0,10){29.5}}
\put(300,30){\line(0,10){9.5}} \put(310,30){\line(0,10){9.5}} \put(320,30){\line(0,10){9.5}} \put(330.5,30){\line(0,10){9.5}} \put(270,20){\rule{10pt
}{10pt}}
\put(270,-10){{\footnotesize $f_{m+3,2,3}(-1)$}}
\thicklines
\put(0,0){\line(10,0){40}} \put(0,0){\line(0,10){40}} \put(0,40){\line(10,0){10}} \put(10,30){\line(10,0){10}} \put(10,30){\line(0,10){10}}
\put(20,20){\line(10,0){10}} \put(20,20){\line(0,10){10}} \put(30,10){\line(10,0){10}} \put(30,10){\line(0,10){10}} \put(40,0){\line(0,10){10}}

\end{picture}

\vspace{10pt}

\section{Preliminaries}

A partition $\lambda$ of a positive integer $n$ is a non-increasing sequence of nonnegative integers $\lambda=(\lambda_1, \lambda_2, \cdots, \lambda_d)$
such that $n=\lambda_1 + \lambda_2 + \cdots + \lambda_d$. A Ferrers diagram of shape $\lambda$ is a left-justified array of $n$ boxes, with row $i$ (from
top to bottom) containing $\lambda_i$ boxes. A standard Young tableau (SYT) of shape $\lambda$ is a labeling by $\{1, 2, \cdots ,n\}$ of the boxes in the Ferrers diagram such that each row and column is increasing (from left to right, and from top to bottom respectively). A skew shape $\lambda/\mu$ is the collection of boxes which belong to $\lambda$ but not $\mu$ when drawn with coinciding top left corners. The number of SYT of shape $\lambda$ is denoted by $f^\lambda$. These types are illustrated as

\begin{picture}(50,50)(-80,0)
\put(0,0){\line(10,0){10}} \put(0,10){\line(10,0){20}} \put(0,20){\line(10,0){40}} \put(0,30){\line(10,0){40}}
\put(0,0){\line(0,10){29.5}} \put(10,0){\line(0,10){29.5}} \put(20,10){\line(0,10){19.5}} \put(30,20){\line(0,10){9.5}} \put(40,20){\line(0,10){9.5}}
\put(0,-10){{\footnotesize $\lambda=(4,2,1)$}}

\put(100,0){\line(10,0){20}} \put(100,10){\line(10,0){40}} \put(110,20){\line(10,0){40}} \put(130,30){\line(10,0){20}}
\put(100,0){\line(0,10){9.5}} \put(110,0){\line(0,10){19.5}} \put(120,0){\line(0,10){19.5}} \put(130,10){\line(0,10){19.5}}
\put(140,10){\line(0,10){19.5}} \put(150,20){\line(0,10){9.5}} \put(90,-10){{\footnotesize $\lambda/\mu=(5,4,2)/(3,1)$}}

\end{picture}

\vspace{10pt}

An excellent source for more on the properties of SYT is \cite{RI}. There are two remarkable formulas on enumerations of $\lambda$ and $\lambda/\mu$.

\vspace{5pt}

{\bf Proposition 2.1} (The Frobenius-Young Formula) \cite{FE,AL} {\it The number of SYT of shape} $\lambda=(\lambda_1,\lambda_2,\cdots,\lambda_d)$ {\it with} $\lambda_1 \geq \lambda_2 \geq \cdots \geq \lambda_d$ {\it is}
$$
f^\lambda = \frac{(\lambda_1 + \cdots + \lambda_d)!}{\prod_{i=1}^d (\lambda_i + d - i)!} \prod_{1 \leq i < j \leq d} (\lambda_i - \lambda_j - i + j). \eqno(2.1)
$$

\vspace{5pt}

{\bf Proposition 2.2} (The Aitken Formula) \cite [Cor. 7.16.3]{RI} {\it The number of skew SYT of shape} $\lambda/\mu=(\lambda_1,\cdots,\lambda_d)/(\mu_1,\cdots,\mu_r)$ {\it such that} $r < d, \mu_i < \lambda_i$,  {\it is}
$$
f^{\lambda/\mu} = n! \; det\left( \frac{1}{(\lambda_i -\mu_j - i + j)!} \right)_{i,j=1}^d, \eqno(2.2)
$$
{\it where} $n=\lambda_1+\cdots+\lambda_d-\mu_1-\cdots-\mu_r$.

The Frobenius-Young Formula is also called hook-length formula \cite{JS}. There are many different elegant proofs of (2.1) in the literature, using generating function \cite{AP}, probabilistics method  \cite{CU}, bijective proofs \cite{DS,JE}, difference methods \cite{PM}, and Schur functions \cite{RI}, respectively. It should be noted that Greene et al. \cite{CU} used discrete probability and all known proofs fall within the scope of discrete mathematics.

We are now going to recall order statistics. Let $X_1, X_2,\cdots,X_n$ be $n$ random variables on some probability space $(\Omega,{\cal F}, P)$, then the corresponding order statistics are obtained by arranging these $n$ $X$'s in nondecreasing order, and are denoted by $X_{1,n}, X_{2,n}, \cdots, X_{n,n}$ such that $X_{1,n} \leq X_{2,n} \leq \cdots \leq X_{n,n}$. For conveniently, we use the order statistics $\xi_{1,n}, \xi_{2,n}, \cdots, \xi_{n,n}$ from i.i.d random variables with uniform distribution on $[0,1]$ in this paper, this turns out the difference between $(\xi_{1,n} \leq \xi_{2,n} \cdots \leq \xi_{n,n})$ and $(\xi_{1,n} < \xi_{2,n} < \cdots < \xi_{n,n})$ is only a measure-zero set.

\vspace{5pt}

{\bf Proposition 2.3} (The probability density function) \cite{HA} {\it Suppose} $\xi_{1,n}, \xi_{2,n}, \cdots, \xi_{n,n}$ {\it are order statistics from i.i.d random variables with uniform distribution on} $[0,1]$, {\it the joint probability density function is}
$$
f(x_1, x_2, \cdots, x_n) = n!, \;\; 0 < x_1 < x_2 < \cdots < x_n < 1. \eqno(2.3)
$$

\section{Probabilistic model of Standard Young tableaux}

The key idea of our model is to establish a kind of continuous background for some discrete structures. Suppose $\lambda_1 \geq \lambda_2 \geq \cdots \geq \lambda_d$, $n=\lambda_1+\cdots+\lambda_d$, there are $d$ groups independent order statistics $(\xi_{1,\lambda_1}, \cdots, \xi_{\lambda_1,\lambda_1})$,  $(\eta_{1,\lambda_2}, \cdots, \eta_{\lambda_2,\lambda_2})$, $\cdots$, $(\tau_{1,\lambda_d}, \cdots, \tau_{\lambda_d,\lambda_d})$. Furthermore, all of these order statistics are assumed to be from i.i.d random variables with uniform distribution on $[0,1]$.

\vspace{5pt}

Fixing $\omega \in \Omega$, it is clear that the following event
$$
A=\left ({\tiny \begin{array}{*{10}l}
  \xi_{1,\lambda_1}(\omega) < \xi_{2,\lambda_1}(\omega) < \cdots  < \xi_{\lambda_2,\lambda_1}(\omega) < \cdots  <\xi_{\lambda_1,\lambda_1}(\omega)  \\
    \;\;\;\;\; \wedge \;\;\;\;\;\;\;\;\;\;\;\;\;\;\; \wedge \;\;\;\;\;\;\;\;\;\; \cdots \;\;\;\;\;\;\;\;\;\; \wedge \\
  \eta_{1,\lambda_2}(\omega) < \eta_{2,\lambda_2}(\omega) < \cdots  < \eta_{\lambda_2,\lambda_2}(\omega)\\
    \;\;\;\;  \cdots \;\;\;\;\;\;\;\; \cdots \;\;\;\;\;\;\;\;\; \cdots \\
    \;\;\;\;\; \wedge \;\;\;\;\;\;\;\; \cdots  \cdots \;\;\;\;\;\;\wedge  \\
   \tau_{1,\lambda_d}(\omega) < \cdots < \tau_{\lambda_d,\lambda_d}(\omega)
\end{array} } \right )  \eqno(3.1)
$$
corresponds to filling the SYT $\lambda=(\lambda_1, \lambda_2,\cdots,\lambda_d)$. It is interesting that there are two completely different ways to evaluating the probability of event $A$ of (3.1).

Let $S_N$ denote the nested simplexes:
$$
S_N=\left ({\tiny \begin{array}{*{10}l}
 0< x_{1,\lambda_1} < x_{2,\lambda_1} < \cdots  < x_{\lambda_2,\lambda_1} < \cdots  <x_{\lambda_1,\lambda_1} < 1  \\
   \;\;\;\;\; \;\;\;\;\; \wedge \;\;\;\;\;\;\;\;\;\; \wedge \;\;\;\;\;\;\;\; \cdots \;\;\;\;\;\;\;\; \wedge \\
 \;\;\;\; \;\;\; y_{1,\lambda_2} < y_{2,\lambda_2} < \cdots  < y_{\lambda_2,\lambda_2} < 1 \\
  \;\;\;\;  \;\;\;\;  \cdots \;\;\;\;\;\;\;\; \cdots \;\;\;\;\;\;\;\;\; \cdots \\
  \;\;\;\;  \;\;\;\;\; \wedge \;\;\;\;\;\;\; \cdots  \cdots \;\;\;\;\;\wedge  \\
  \;\;\;\; \;\;\; z_{1,\lambda_d} < \cdots < z_{\lambda_d,\lambda_d} < 1
\end{array} } \right ),  \eqno(3.2)
$$
(2.3) follows that
$$
{\bf P}(A)=\prod_{i=1}^d \lambda_i! \int\int_{S_N} dx_{1,\lambda_1}\cdots dz_{\lambda_d,\lambda_d}.
$$

On the other hand, the property of i.i.d makes sure that all the outcomes of (3.1) are equally likely. Based on this discrete structure there is
$$
{\bf P}(A) = \frac{f^\lambda}{{n \choose \lambda_1, \lambda_2, \cdots, \lambda_d}}.
$$

\vspace{5pt}

Therefore, we have the following surprising result

\vspace{5pt}

{\bf Proposition 3.1} {\it The number of SYT $f^\lambda$ is}
$$
f^\lambda = {n \choose \lambda_1, \lambda_2, \cdots, \lambda_d} {\bf P} (A) = n! \int \int_{S_N} dx_{1,\lambda_1} \cdots dz_{\lambda_d,\lambda_d}. \eqno(3.3)
$$

It should be noted that the order statistics model is suitable for skew type $\lambda/\mu$ and other types (shifted SYT, truncated shape, etc.). In this section we give three examples.

\vspace{5pt}

{\bf Example 3.1}. It is well-known $f^{(4,2,1)} = \frac{7!}{6 \cdot 4 \cdot 2 \cdot 1 \cdot 3 \cdot 1 \cdot 1} = 35$. From (3.3) one may derive it as follows
\begin{eqnarray*}
f^{(4,2,1)} & = & {7 \choose 4 \;\; 2 \;\; 1} {\bf P} \left ( {\tiny \begin{array}{ccccccc}
\xi_{1,4} & < & \xi_{2,4} & < & \xi_{3,4} & < & \xi_{4,4} \\
\wedge    &   & \wedge    &   &   &   &   \\
\eta_{1,2}& < & \eta_{2,2} &   &   &   &   \\
\wedge    &   &     &   &   &   &   \\
\tau_{1,1}&   &   &   &   &   &
                                                      \end{array}}
 \right ) \\
& = & 7! \mathop {\int {\int {...\int {} } } }\limits_{\tiny \begin{array}{*{10}l}
   0 < x_1 < x_2 <  x_3 < x_4 < 1  \\
    \;\;\;\;\;\;\;\; \wedge \;\;\;\;\; \wedge   \\
    \;\;\;\;\;\;\; y_1 < y_2 < 1 \\
    \;\;\;\;\;\;\;\; \wedge   \\
   \;\;\;\;\;\;\;z < 1
\end{array} } \; dx_1 \cdots dx_4 dy_1 dy_2 dz \\
& = & 7! \mathop {\int {\int {...\int {} } } }\limits_{\tiny \begin{array}{*{10}l}
   0 < x_1 < x_2   \\
    \;\;\;\;\;\;\;\; \wedge \;\;\;\;\; \wedge   \\
    \;\;\;\;\;\;\; y_1 < y_2 < 1 \\
\end{array} } \; (1-y_1)\frac{(1-x_2)^2}{2!} \; dx_1 \cdots dy_2 = \frac{7!}{3!} \int_0^1 \frac{7}{24} t^6 dt = 35.
\end{eqnarray*}

\vspace{5pt}

{\bf Example 3.2}. The number of skew type $(5,4,2)/(3,1)$ may be evaluated as follows
\begin{eqnarray*}
f^{(5,4,2)/(3,1)} & = & {7 \choose 2 \;\; 3 \;\; 2} {\bf P} \left ( {\tiny \begin{array}{ccccccccc}
          &   &             &   &   &   & \xi_{1,2} & < & \xi_{2,2} \\
          &   &  &  &  &                    & \wedge  &   &   \\
          &   &\eta_{1,3}   & < & \eta_{2,3} & <  & \eta_{3,3}  &   &   \\
        &    & \wedge      &     &   &   &   &   \\
\tau_{1,1}& < & \tau_{1,2}  &   &   &   &           &   &
                                                      \end{array}}
 \right ) \\
& = & 7! \mathop {\int {\int {...\int {} } } }\limits_{\tiny \begin{array}{*{10}l}
   \;\;\;\;\;\;\;\;\; \;\;\;\;\;\;\;\;\;\;\;\;\;\;\;0 < x_1 < x_2 <   1  \\
   \;\;\;\; \;\;\;\;\;\;\;\;\;\;\;\;\;\;\; \;\;\;\;\;\;\; \;\;\;\;\; \wedge   \\
   \; \;\;\;\;\;\;\; 0 < y_1 < y_2 < y_3 < 1 \\
    \;\;\;\;\;\;\;\;\;\;\;\;\;\;\; \wedge   \\
  0 < z_1 < z_2 < 1
\end{array} } \; dx_1 dx_2 dy_1 dy_2 dy_3 dz_1 dz_2 \\
& = & 7! \mathop {\int {\int {\int {} } } }\limits_{\tiny  0 < y_1 < y_2 < y_3 < 1   } \; \frac{(1-y_1^2)}{2!} \times \frac{1-(1-y_3)^2}{2!} \; dy_1 dy_2 dy_3 = 169.
\end{eqnarray*}

\vspace{5pt}

{\bf Example 3.3}. Similarly, the following multiple integral
$$
7! \mathop {\int {\int {...\int {} } } }\limits_{\tiny \begin{array}{*{10}l}
   0 < x_1 < x_2 <  x_3 < x_4 < 1  \\
   \;\;\;\;\;\;\;\; \;\;\;\;\;\;\;\; \wedge \;\;\;\;\; \wedge   \\
   \;\;\;\;\;\;\;\; \;\;\;\;\;\;\; y_1 < y_2  \\
  \;\;\;\;\;\; \; \;\;\;\;\;\;\;\; \; \;\;\;\;\;\;\;\; \wedge   \\
  \;\;\;\;\;\; \;\;\;\;\;\;\;\;\; \;\;\;\;\;\;\;\;z < 1
\end{array} } \; dx_1 \cdots dx_4 dy_1 dy_2 dz = 7
$$
gives the known number $g^{(4,2,1)}=7$ of the shifted SYT \hspace{2cm}.
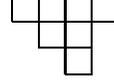
\begin{picture}(0,0)(60,0)
\put(20,0){\line(10,0){10}} \put(10,10){\line(10,0){20}} \put(0,20){\line(10,0){40}} \put(0,30){\line(10,0){40}}
\put(0,20){\line(0,10){9.5}} \put(10,10){\line(0,10){19.5}} \put(20,0){\line(0,10){29.5}} \put(30,0){\line(0,10){29.5}} \put(40,20){\line(0,10){9.5}}
\end{picture}

\section{Main results}

In the following lemma, which will be useful in the main result of this paper, we investigate the integral of an order statistics.

\vspace{5pt}

{\bf Lemma 4.1.} (Determinant of SYT) {\it For } $\lambda_1 \geq \lambda_2 \geq \cdots \geq \lambda_d > 0, \; 0 < t_1 \leq t_2 \leq \cdots \leq t_d \leq 1 $, {\it we have the following multiple integral formula}
$$
J_{\lambda}(t_1, \cdots, t_d)=  \mathop {\int {\int {...\int {} } } }\limits_{\tiny \begin{array}{*{10}l}
   0<x_{1,1} < x_{1,2} < \cdots <x_{1,\lambda_1} < t_1  \\
    \;\;\;\;\;\;\;\; \wedge \;\;\;\;\;\;\;\; \wedge \;\; \cdots   \\
    \;\;\;\;\;\;\;\;\;\;\; \cdots \cdots \cdots \cdots\\
    \;\;\;\;\;\;\;\; \wedge \;\;\; \cdots  \cdots \;\;\;\; \wedge  \\
   \;\;\;\;\;\;x_{d,1} < \cdots < x_{d,\lambda_d} < t_d
\end{array} } \; dx_{1,1} \cdots dx_{1,\lambda_1} \cdots dx_{d,1} \cdots dx_{d,\lambda_d}
$$
$$
\hspace{-50pt}=det\left( \frac{(\lambda_i+d-i)_{d-j} \; t_i^{\lambda_i + j - i}}{(\lambda_i + d - i)!} \right)_{i,j=1}^d. \eqno(4.1)
$$

\vspace{5pt}

{\bf Proof.} It suffices to show that for $0 < y_1 \leq y_2 \leq \cdots \leq y_d <1$,
$$
F_d(y_1,\cdots,y_d)= \mathop {\int {\int {...\int {} } } }\limits_{\tiny \begin{array}{*{10}l}
   0<x_1 < x_2 < \cdots <x_d  \\
   \;\;\;\;\;\;\; \wedge \;\;\;\;\;\; \wedge \;\;\;\; \cdots \;\;\;\; \wedge  \\
   \;\;\;\;\;\;\;y_1 \;\;\;\;\; y_2 \;\;\; \cdots \;\;\;\; y_d
\end{array} } \; dx_1 \cdots dx_d = det \left( \frac{(j)_{i-1} \; y_i^{j-i+1}}{j!} \right )_{i,j=1}^d
$$
$$
= \left | \begin{array}{ccccccc}
   y_1 & \frac{y_1^2}{2!} & \frac{y_1^3}{3!} & \cdots &          &         & \frac{y_1^d}{d!} \\
   1   & y_2              & \frac{y_2^2}{2!} & \cdots &          &         & \frac{y_2^{d-1}}{(d-1)!} \\
   0   & 1                & y_3              & \cdots &          &         & \frac{y_2^{d-2}}{(d-2)!}  \\
       & \cdots           &                  &        & \cdots   &         &  \cdots   \\
       &                  &                  &        &  1       & y_{d-1} & \frac{y_{d-1}^2}{2!} \\
       &   \cdots         &                  &        &  0       & 1       & y_d
 \end{array} \right |.   \eqno(4.2)
$$
We observe that
$$
F_d(y_1,\cdots,y_d)= \mathop {\int {\int {...\int {} } } }\limits_{\tiny \begin{array}{*{10}l}
   0<x_1 < x_2 < \cdots <x_{d-1}  \\
   \;\;\;\;\;\;\; \wedge \;\;\;\;\;\; \wedge \;\;\;\; \cdots \;\;\;\; \wedge  \\
   \;\;\;\;\;\;\;y_1 \;\;\;\;\; y_2 \;\;\; \cdots \;\;\;\; y_{d-1}
\end{array} } \; (y_d - x_{d-1}) dx_1 \cdots dx_{d-1}
$$
\begin{eqnarray*}
& = & y_d \cdot F_{d-1}(y_1, \cdots, y_{d-1}) - \mathop {\int {\int {...\int {} } } }\limits_{\tiny \begin{array}{*{10}l}
   0<x_1 < x_2 < \cdots <x_{d-2}  \\
   \;\;\;\;\;\;\; \wedge \;\;\;\;\;\; \wedge \;\;\;\; \cdots \;\;\;\; \wedge  \\
   \;\;\;\;\;\;\;y_1 \;\;\;\;\; y_2 \;\;\; \cdots \;\;\;\; y_{d-2}
\end{array} } \; (\frac{y_{d-1}^2}{2!} - \frac{x_{d-2}^2}{2!}) dx_1 \cdots dx_{d-2} \\
& = & \cdots \cdots \\
& = & y_dF_{d-1} - \frac{y_{d-1}^2}{2!}F_{d-2} + \frac{y_{d-2}^3}{3!} F_{d-3} - \cdots + (-1)^d \frac{y_2^{d-1}}{(d-1)!} y_1 + (-1)^{d+1}\frac{y_1^d}{d!},
\end{eqnarray*}
which is the determinant (4.2) by expanding along the last column.

\vspace{5pt}

Notice that the $i$-th row of determinant (4.2) is related only to $y_i$, we can add the $(d-1)$-th row to the last row after integrating with respect to $y_d$, and then add the $(d-2)$-th row to the $(d-1)$-th row after integrating with respect to $y_{d-1}$, etc. So it follows
$$
\mathop {\int {\int {...\int {} } } }\limits_{\tiny \begin{array}{*{10}l}
   0<y_1 < y_2 < \cdots <y_d  \\
   \;\;\;\;\;\;\; \wedge \;\;\;\;\;\; \wedge \;\;\;\; \cdots \;\;\;\; \wedge  \\
   \;\;\;\;\;\;\;z_1 \;\;\;\;\; z_2 \;\;\; \cdots \;\;\;\; z_d
\end{array} } \; F_d(y_1,\cdots,y_d) \; dy_1 \cdots dy_d
$$
$$
= \left | \begin{array}{ccccccc}
\frac{z_1^2}{2!} & \frac{z_1^3}{3!}   & \frac{z_1^4}{4!} & \cdots &          &                      & \frac{z_1^{d+1}}{(d+1)!} \\
   z_2           & \frac{z_2^2}{2!}   & \frac{z_2^3}{3!} & \cdots &          &                      & \frac{z_2^d}{d!} \\
   1             & z_3                & \frac{z_3^2}{2!} & \cdots &          &                      & \frac{z_3^{d-1}}{(d-1)!}  \\
                 & \cdots             &                  &        & \cdots   &                      &  \cdots   \\
                 &                    &                  &        &  z_{d-1} & \frac{z_{d-1}^2}{2!} & \frac{z_{d-1}^3}{3!} \\
                 &   \cdots           &                  &        &  1       & z_d                  & \frac{z_d^2}{2!}
 \end{array} \right |.
$$
Therefore (4.1) follows from the fact that $J_{\lambda}(t_1, \cdots, t_d)$ is equal to
$$
\mathop {\int {\int {...\int {} } } }\limits_{\tiny \begin{array}{*{10}l}
   0<x_{1,2} < x_{1,3} < \cdots <x_{1,\lambda_1} < t_1  \\
    \;\;\;\;\;\;\;\; \wedge \;\;\;\;\;\;\;\; \wedge \;\; \cdots   \\
    \;\;\;\;\;\;\;\;\;\;\; \cdots \cdots \cdots \cdots\\
    \;\;\;\;\;\;\;\; \wedge \;\;\; \cdots  \cdots \;\;\;\; \wedge  \\
   \;\;\;\;\;\;x_{d,2} < \cdots < x_{d,\lambda_d} < t_d
\end{array} } F_d(x_{1,2}, \cdots, x_{d,2}) \; dx_{1,2} \cdots dx_{1,\lambda_1} \cdots dx_{d,2} \cdots dx_{d,\lambda_d}
$$
which order of integration is $x_{d,2}, \cdots, x_{1,2}, x_{d,3}, \cdots, x_{1,3}, \cdots, x_{1,\lambda_1}$. \,\,\,\,\, $\Box$

\vspace{5pt}

{\bf Corollary 4.2}. {\it The Frobennius-Young formula (2.1) and Aitken formula (2.2) hold true}.

\vspace{5pt}

{\bf Proof}. Write $x_i=\lambda_i+d-i, 1 \leq i \leq d$, then $x_1 > x_2 > \cdots > x_d$. From (4.1) we have
$$
J_{\lambda}(1, \cdots, 1) = \frac{ det\left((x_i)_{d-j} \right)_{i,j=1}^d}{\prod_{i=1}^d x_i!}.
$$
It is clear that $det((x_i)_{d-j})_{i,j=1}^d$ is Vandermonde determinant from the fact that
$$
(x)_r = \sum_{k=0}^r s(r,k) x^k,
$$
here $s(r,k)$ is Stirling numbers of the first kind. Hence Frobennius-Young formula (2.1) follows from (3.3) and (4.1) immediately.

By making minor changes in the proof of lemma 4.1, we can derive Aitken formula (2.2) directly. Note that the left $\mu_r$ columns of the skew SYT $(\lambda_1,\cdots,\lambda_d)/(\mu_1,\cdots,\mu_r)$ construct the SYT of shape $(\mu_r,\alpha_2,\cdots,\alpha_{d-r})$, here $\alpha_j=min(\mu_r, \lambda_{r+j})$. From (4.1) the related determinant of $(\mu_r,\alpha_2,\cdots,\alpha_{d-r})$ is
$$
B = \left | \begin{array}{cccc}
\frac{x_p^{\mu_r}}{\mu_r!} & \frac{x_p^{\mu_r+1}}{(\mu_r+1)!}       & \cdots  & \frac{x_p^{\mu_r+d-r-1}}{(\mu_r+d-r-1)!} \\
       \cdot               & \frac{x_{p+1}^{\alpha_2}}{\alpha_2!}   & \cdots   & \cdots  \\
       \cdot               &                                         &         &          \\
        \cdot               & &         &   \cdots
 \end{array} \right |_{(d-r) \times (d-r)}.
$$
To complete the similar integration, it is necessary to extend this determinant to be
{\footnotesize $$
 \left | \begin{array}{cccccccccc}
1   & x_1 & \cdots & \frac{x_1^{p-4}}{(p-4)!} & \frac{x_1^{p-3}}{(p-3)!}  & \frac{x_1^{p-2}}{(p-2)!}  &  | & \frac{x_1^{\mu_r+p-1}}{(\mu_r+p-1)!} &
 \frac{x_1^{\mu_r+p}}{(\mu_r+p)!}  & \cdots \\
0   & 1   &  \cdots &    &  \cdots  &     &  |  & \cdots &     &      \\
0   & 0   &    &  1  &  x_{p-3} &  \frac{x_{p-3}^2}{2!} & | & \frac{x_{p-3}^{\mu_r+3}}{(\mu_r+3)!} & \frac{x_{p-3}^{\mu_r+4}}{(\mu_r+4)!} & \cdots \\
0   & \cdots & &  0  &  1 & x_{p-2} & | & \frac{x_{p-2}^{\mu_r+2}}{(\mu_r+2)!} & \frac{x_{p-2}^{\mu_r+3}}{(\mu_r+3)!} & \cdots \\
0   & \cdots & &  0  &  0 &  1      & | & \frac{x_{p-1}^{\mu_r+1}}{(\mu_r+1)!} & \frac{x_{p-1}^{\mu_r+2}}{(\mu_r+2)!} & \cdots \\
 -  & -  &  - & - & - & - & -  & - & - &  - \\
    &    &    &   &   &   & | &   &   &  \\
    &    &   & 0  &   &   & | &  B  & &    \\
    &    &    &   &   &   & | &   &   &
 \end{array} \right |
$$}
in the case of
\begin{picture}(0,50)(-20,-10)
\put(-10,10){\line(10,0){100}} \put(-10,-10){\line(0,10){20}} \put(70,0){\line(0,10){30}}
\put(100,40){\line(10,0){30}} \put(70,30){\line(10,0){30}} \put(100,30){\line(0,10){10}}
\put(75,2){{\footnotesize $x_p$}}   \put(75,21){{\footnotesize $x_1$}} \put(78,15){{\footnotesize $\cdot$}} \put(78,11){{\footnotesize $\cdot$}}
\put(-5,-5){{\footnotesize $(\mu_r,\alpha_2,\cdots,\alpha_{d-r})$}}
\end{picture} \hspace{5cm}.

By a similar argument, we obtain Aitken formula (2.2). \,\,\,\,\,\, $\Box$

{\bf Remark 4.1}. The determinant (4.1), which diagonal entries are $t_i^{\lambda_i}/\lambda_i!$,
$$
 \left | \begin{array}{ccccccc}
\frac{t_1^{\lambda_1}}{\lambda_1!}         & \frac{t_1^{\lambda_1+1}}{(\lambda_1+1)!} & \frac{t_1^{\lambda_1+2}}{(\lambda_1+2)!} & \cdots &     & \cdots  & \frac{t_1^{\lambda_1+d-1}}{(\lambda_1+d-1)!} \\
\frac{t_2^{\lambda_2-1}}{(\lambda_2-1)!}   & \frac{t_2^{\lambda_2}}{\lambda_2!}       & \frac{t_2^{\lambda_2+1}}{(\lambda_2+1)!} & \cdots &     & \cdots & \frac{t_2^{\lambda_2+d-2}}{(\lambda_2+d-2)!} \\
\frac{t_3^{\lambda_3-2}}{(\lambda_3-2)!}   & \frac{t_3^{\lambda_3-1}}{(\lambda_3-1)!} & \frac{t_3^{\lambda_3}}{\lambda_3!}       & \cdots &     &  \cdots & \frac{t_3^{\lambda_3+d-3}}{(\lambda_3+d-3)!}  \\
       & \cdots   &          &        & \cdots   &         &  \cdots   \\
0 & 1 & t_{d-1} & \frac{t_{d-1}^2}{2!} & \cdots & \frac{t_{d-1}^{\lambda_{d-1}}}{(\lambda_{d-1})!} & \frac{t_{d-1}^{\lambda_{d-1}+1}}{(\lambda_{d-1}+1)!} \\
\cdots  &   0 &  1  &  t_d   &  \cdots       & \frac{t_d^{\lambda_d-1}}{(\lambda_d-1)!}            & \frac{t_d^{\lambda_d}}{\lambda_d!}
 \end{array} \right |
$$
is constructed by placing $\frac{t_i^{\lambda_i+d-i}}{(\lambda_d+d-i)!}(1 \leq i \leq d)$ in the last column, the first derivative of each term in the $(d-1)$-th column, and so on through the $(d-1)$st derivative.

%One may construct (4.1) following the rules:
%(1). {\bf Row Rule}. {\it The $(i,i-k)$-th entry is ${\bf D}^k (t_i^{\lambda_i}/\lambda_i!)$, {\bf D} is differential operator}.
%(2). {\bf Column Rule}. {\it The $(i-k,i)$-th entry is $t_{i-k}^{\lambda_{i-k}+k}/(\lambda_{i-k}+k)!)$}.

\vspace{5pt}

{\bf Remark 4.2}. The determinant of skew SYT of shape $\lambda/\mu$ is
$$
J_{\lambda/\mu}(t_1,\cdots,t_d)=det\left( \frac{(\lambda_i+d-i)_{d+\mu_j-j} \; t_i^{\lambda_i -\mu_j + j -i}}{(\lambda_i + d - i)!} \right)_{i,j=1}^d. \eqno(4.3)
$$

\vspace{5pt}

Moreover, we give the following example to explain (4.3).

\vspace{5pt}

{\bf Example 4.1}. From (2.2), the number of skew SYT of shape $(6,5,4,4,1)/(3,2,2)$ is
$$
f^{(6,5,4,4,1)/(3,2,2)} = 13! \left | \begin{array}{ccccc}
\frac{1}{3!} & \frac{1}{5!} & \frac{1}{6!} & \frac{1}{9!} & \frac{1}{10!} \\
   1         & \frac{1}{3!} & \frac{1}{4!} & \frac{1}{7!} & \frac{1}{8!} \\
   0         &      1       & \frac{1}{2!} & \frac{1}{5!} & \frac{1}{6!} \\
   0         &      1       &      1       & \frac{1}{4!} & \frac{1}{5!} \\
   0         &      0       &      0       &       1      &     1
 \end{array} \right |.
$$
In fact, $J_{(6,5,4,4,1)/(3,2,2)}(t_1, \cdots,t_5)$ is equal to
$$
\mathop {\int {\int {...\int {} } } }\limits_{\tiny \begin{array}{*{10}l}
 \;\; \;\;\;\;\;\; 0  < y_0 < z_0  <    s < t_1  \\
 \;\;  \;\;\;\;\;\; \;\;\;\;\;\;\; \wedge \;\;\;\;\;\; \wedge    \\
 0 < x_1 < y_1 < z_1 < t_2 \\
    \;\;\;\;\;\;\; \wedge \;\;\; \;\;\; \wedge  \\
   \;\;\;\;\;\;x_2 < y_2 < t_3 \\
    \;\;\;\;\;\;\; \wedge \;\;\; \;\;\; \wedge  \\
   \;\;\;\;\;\;x_3 < y_3 < t_4
\end{array} } \left | \begin{array}{cc}
\frac{x_3^2}{2!} & \frac{x_3^3}{3!}  \\
   1             &  t_5
 \end{array} \right | \; dx_1x_2x_3 dy_0 \cdots y_3 dz_0 z_1 ds.
$$
Clearly,
$$
\left | \begin{array}{cc}
\frac{x_3^2}{2!} & \frac{x_3^3}{3!}  \\
   1             &  t_5
 \end{array} \right | = \left | \begin{array}{cccc}
1  &  x_1  & \frac{x_1^4}{4!}  &  \frac{x_1^5}{5!} \\
0  &    1  & \frac{x_2^3}{3!}  &  \frac{x_2^4}{4!} \\
0  &    0  & \frac{x_3^2}{2!} & \frac{x_3^3}{3!}  \\
0  &    0  &   1             &  t_5
 \end{array} \right | ,
$$
and
$$
\left | \begin{array}{cccc}
y_1  & \frac{y_1^2}{2!}  & \frac{y_1^5}{5!}  &  \frac{y_1^6}{6!} \\
1  &    y_2  & \frac{y_2^4}{4!}  &  \frac{y_2^5}{5!} \\
0  &    1  & \frac{y_3^3}{3!} & \frac{y_3^4}{4!}  \\
0  &    0  &   1             &  t_5
 \end{array} \right | = \left | \begin{array}{ccccc}
 1 & \frac{y_0^2}{2!} & \frac{y_0^3}{3!} & \frac{y_0^6}{6!} & \frac{y_0^7}{7!}  \\
 0 & y_1  & \frac{y_1^2}{2!}  & \frac{y_1^5}{5!}  &  \frac{y_1^6}{6!} \\
 0 & 1  &    y_2  & \frac{y_2^4}{4!}  &  \frac{y_2^5}{5!} \\
 0 & 0  &    1  & \frac{y_3^3}{3!} & \frac{y_3^4}{4!}  \\
 0 & 0  &    0  &   1             &  t_5
 \end{array} \right | .
$$
Therefore, the number of skew SYT of shape $(6,5,4,4,1)/(3,2,2)$ follows from
$$
J_{(6,5,4,4,1)/(3,2,2)}(t_1, \cdots,t_5) =  \left | \begin{array}{ccccc}
\frac{t_1^3}{3!} & \frac{t_1^5}{5!} & \frac{t_1^6}{6!} & \frac{t_1^9}{9!} & \frac{t_1^{10}}{10!} \\
   t_2         & \frac{t_2^3}{3!} & \frac{t_2^4}{4!} & \frac{t_2^7}{7!} & \frac{t_2^8}{8!} \\
   0         &      t_3       & \frac{t_3^2}{2!} & \frac{t_3^5}{5!} & \frac{t_3^6}{6!} \\
   0         &      1       &      t_4       & \frac{t_4^4}{4!} & \frac{t_4^5}{5!} \\
   0         &      0       &      0       &       1      &     t_5
 \end{array} \right |.
$$

\section{Several special truncated shapes}

Using complex analysis, Adin et al. \cite{RO} and Panova \cite{GR} discussed the number of rectangle truncated by a staircase. In this section we use (3.3) and lemma 4.1 to research SYT of truncated shapes.

Write $f_{\alpha_1, \alpha_2, \cdots, \alpha_k}$ to be the number of a left-justified SYT of truncated shape $\lambda_T = (\alpha_1, \alpha_2, \cdots, \alpha_k)$, with $\alpha_i$ boxes in the $i$-th row. For examples, $f_{m+1, m+2, \cdots, m+k}$ is the number of a rectangle $(m+k)^k$ removed the right $k-i$ boxes in the $i$-th row ($1 \leq i \leq k-1$), $f_{m+k,m-1,m}$ is the number of SYT $\lambda=(m+k,m,m)$ truncated by the last box in the second row. The following Theorem 5.1 is a specialization of Theorem 2 in \cite{GR}.

\vspace{5pt}

{\bf Theorem 5.1}. {\it The number of SYT of truncated shape $\lambda_T=(m+1, m+2, \cdots, m+k)$ is}
$$
f_{m+1,m+2,\cdots,m+k}=\left ( mk+ {k+1 \choose 2} \right )! \prod_{i=0}^{k-1} \frac{\Gamma(m+1+\frac{i}{2}) \Gamma(1+\frac{i}{2}) \Gamma(\frac{1+i}{2})}{\Gamma(m+1+\frac{i+k+1}{2}) \Gamma(m+1+i) \Gamma(\frac{1}{2})}. \eqno(5.1)
$$

{\bf Proof}. From (3.3), Lemma 4.1 yields that
$$
f_{m+1,m+2,\cdots,m+k}=\left ( mk+ {k+1 \choose 2} \right )! \; J_{(m+1,\cdots,m+k)}(1,\cdots,1)
$$
where
$$
J_{(m+1,\cdots,m+k)}(1,\cdots,1)=\mathop {\int{...\int}} \limits_{\tiny 0< x_1 < x_2 < \cdots < x_k < 1} \left | \begin{array}{cccc}
\frac{x_1^m}{m!}  & \frac{x_1^{m+1}}{(m+1)!}  & \cdots &  \frac{x_1^{m+k-1}}{(m+k-1)!} \\
\frac{x_2^m}{m!}  & \frac{x_2^{m+1}}{(m+1)!}  & \cdots &  \frac{x_2^{m+k-1}}{(m+k-1)!}  \\
&  \cdots   & \cdots &  \\
\frac{x_k^m}{m!}  & \frac{x_k^{m+1}}{(m+1)!}  & \cdots &  \frac{x_k^{m+k-1}}{(m+k-1)!}
 \end{array} \right | dx_1 \cdots dx_k
$$
$$
= \mathop {\int{...\int}} \limits_{\tiny 0< x_1 < x_2 < \cdots < x_k < 1} \prod_{i=1}^k \frac{x_i^m}{(m+i-1)!} \; \prod_{1 \leq i < j \leq k} (x_j - x_i) dx_1 \cdots dx_k.
$$

Immediately, (5.1) follows from (2.3) and the following Selberg's integral formula \cite{PE}
$$
\int_0^1 \cdots \int_0^1 \prod_{i=1}^k t_i^{\alpha-1}(1-t_i)^{\beta-1} \; \prod_{1 \leq i < j \leq k} |t_i-t_j|^{2\gamma} dt_1\cdots dt_k
$$
$$
= k! \prod_{i=0}^{k-1} \frac{\Gamma(\alpha+ i \gamma) \Gamma(\beta+i \gamma) \Gamma(\gamma+i \gamma)}{\Gamma(\alpha+\beta+(i+k-1)\gamma) \Gamma(\gamma)},  \eqno(5.2)
$$
$$
\hspace{50pt} {\rm Re}(\alpha) > 0, {\rm Re}(\beta) > 0, {\rm Re}(\gamma) > -min \{\frac{1}{k}, \frac{{\rm Re}(\alpha)}{k-1}, \frac{{\rm Re}(\beta)}{k-1} \}
$$
with $\alpha=m+1, \beta=1, \gamma=\frac{1}{2}$. \,\,\,\,\,\, $\Box$

\vspace{10pt}

{\bf Remark 5.1}. Furthermore, Selberg's integral gives the number of partially standard Young tableau (PSYT) $(m+1,\cdots, m+k) \subset (m+r+1,\cdots, m+r+k)$. Suppose that $A=(a_1,\cdots,a_d), B=(b_1,\cdots, b_d), a_i \leq b_i, n=b_1+\cdots +b_d$. For a left-justified array of $n$ boxes, with row $i$ (from top to bottom) containing $b_i$ boxes, a PSYT of shape $A \subset B$ is a labeling by $\{1, \cdots, n \}$ of the boxes such that each row is increasing from left to right, but only columns in $A$ is increasing from top to bottom.

A PSYT of shape $A \subset B$ is constructed by adding $b_i-a_i$ boxes to the end of row $i$ $(1 \leq i \leq d
)$ for SYT of (truncated) shape $A$ such that each row is increasing. Therefore, for $m, r \geq 0, k \geq 1$, the number of PSYT of shape $(m+1,\cdots, m+k) \subset (m+r+1,\cdots, m+r+k)$ is
$$
\left ( (m+r)k+ {k+1 \choose 2} \right )! \prod_{i=0}^{k-1} \frac{\Gamma(m+1+\frac{i}{2}) \Gamma(r+1+\frac{i}{2}) \Gamma(\frac{1+i}{2})}{\Gamma(m+r+1+\frac{i+k+1}{2}) \Gamma(m+1+i) \Gamma(r+1) \Gamma(\frac{1}{2})}. \eqno(5.3)
$$
For example, $f_{(2,3,4)}=12,\; f_{(1,2,3) \subset (2,3,4)} = 72$.
\begin{picture}(50,50)(145,10)

\put(150,10){\line(10,0){40}} \put(150,20){\line(10,0){40}} \put(150,30){\line(10,0){30}} \put(150,40){\line(10,0){20}}
\put(150,10){\line(0,10){29.5}} \put(160,10){\line(0,10){29.5}} \put(170,10){\line(0,10){29.5}} \put(180,10){\line(0,10){19.5}}
\put(190.5,10){\line(0,10){10}}
\put(153,31){{\footnotesize $1$}} \put(163,31){{\footnotesize $2$}} \put(153,21){{\footnotesize $3$}} \put(163,21){{\footnotesize $5$}}
\put(173,21){{\footnotesize $9$}} \put(153,11){{\footnotesize $4$}} \put(163,11){{\footnotesize $6$}} \put(173,11){{\footnotesize $7$}}
\put(183,11){{\footnotesize $8$}}
\thicklines
\put(150,40){\line(10,0){10}} \put(160,30){\line(10,0){10}} \put(170,20){\line(10,0){10}}
\put(160,30){\line(0,10){10}} \put(170,20){\line(0,10){10}} \put(180,10){\line(0,10){10}}
\put(150,10){\line(10,0){30}} \put(150,10){\line(0,10){30}}

\end{picture}
is a filling of PSYT of shape $(1,2,3) \subset (2,3,4)$, not a filling of SYT of shape $(2,3,4)$.

\vspace{10pt}

{\bf Theorem 5.2}.

{\it (1). For $m \geq 1, 0 \leq k \leq m$, the number of SYT of shape $(m+1,m+1,m-k)$ with NE corner box removed, is }
$$
f_{m,m+1,m-k} = {3m-k+1 \choose m \;\; m+1 \;\; m-k} \frac{(k+2)^2(2m+1) + (m+2)}{(m+1)(m+2)(2m+1)(2m+3)}. \eqno(5.4)
$$

\vspace{5pt}

{\it (2). For $m \geq 2, k \geq 0$, the number of SYT of shape $(m+k,m,m)$ with last box in the second row removed, is}
{\scriptsize$$
f_{m+k,m-1,m} =  \frac{1}{(2m-1)(2m-3)} \left [ \frac{(k+2)^2(2m-3)+m}{(m+k+2)(m+k+1)} {3m+k-1 \choose m+k \; m-1 \; m}
- \frac{3}{3m-1}{3m-1 \choose m+1 \; m-2 \; m} \right ].  \eqno(5.5)
$$}

{\it (3). For $m \geq 0$, the number of SYT of shape $(m+3,3,3)$ with middle box in the second row removed, is}
$$
f_{m+3,2,3(-1)} = \frac{m+5}{10} {m+2 \choose 2} {m+9 \choose 2}. \eqno(5.6)
$$

\vspace{5pt}

{\bf Proof}. (1). From (3.3) and Lemma 1 we have
$$
f_{m,m+1,m-k}=(3m-k+1)! \mathop {\int{\int}} \limits_{\tiny 0< x < y <  1} \left | \begin{array}{ccc}
\frac{x^{m-1}}{(m-1)!}  & \frac{x^m}{m!}  &   \frac{x^{m+1}}{(m+1)!} \\
\frac{y^{m-1}}{(m-1)!}  & \frac{y^m}{m!}  &   \frac{y^{m+1}}{(m+1)!} \\
\frac{(m-k)_2}{(m-k)!}  & \frac{m-k}{(m-k)!}  &  \frac{1}{(m-k)!}
 \end{array} \right | dx dy,
$$
here the integral of determinant is equal to
$$
 \left | \begin{array}{cc}
\frac{{2m \choose m}}{(2m+1)!}  & \frac{{2m+1 \choose m}}{(2m+2)!}   \\
\frac{m-k}{(m-k)!}  &  \frac{1}{(m-k)!}
 \end{array} \right | - \left | \begin{array}{cc}
\frac{{2m \choose m+1}}{(2m+1)!}  & \frac{{2m+2 \choose m+1}}{(2m+3)!}   \\
\frac{(m-k)_2}{(m-k)!}  &  \frac{1}{(m-k)!}
 \end{array} \right | + \left | \begin{array}{cc}
\frac{{2m+1 \choose m+2}}{(2m+2)!}  & \frac{{2m+2 \choose m+2}}{(2m+3)!}   \\
\frac{(m-k)_2}{(m-k)!}  &  \frac{m-k}{(m-k)!}
 \end{array} \right |
$$
$$
=\frac{1}{m!(m+1)!(m-k)!} \left[\frac{1}{2m+1} - \frac{(m-k)(m+k+4)}{(m+1)(m+2)(2m+3)} \right ].
$$

\vspace{5pt}

(2). The associated multiple integral is
$$
\mathop {\int {\int {...\int {} } } }\limits_{\tiny \begin{array}{*{10}l}
   0 < x_0 < \; x   \\
    \;\;\;\;\;\;\;\; \wedge \;\;\;\;\;\; \wedge   \\
    \;\;\;\;\;\;\; y_0 \; < y < 1 \\
\end{array} } \; \left | \begin{array}{ccc}
\frac{x_0^{m-2}}{(m-2)!}  & \frac{x_0^{m-1}}{(m-1)!}  &   \frac{x_0^m}{m!} \\
\frac{y_0^{m-3}}{(m-3)!}  & \frac{y_0^{m-2}}{(m-2)!}  &   \frac{y_0^{m-1}}{(m-1)!} \\
\frac{y^{m-3}}{(m-3)!}  & \frac{y^{m-2}}{(m-2)!}  &   \frac{y^{m-1}}{(m-1)!}
 \end{array} \right |\frac{(1-x)^k}{k!} \; dx_0 dx dy_0 dy
$$
$$
=\mathop {\int{\int}} \limits_{\tiny 0< x < y <  1} \left | \begin{array}{ccc}
\frac{x^{m-1}}{(m-1)!}  & \frac{x^m}{m!}  &   \frac{x^{m+1}}{(m+1)!} \\
\frac{y^{m-2}}{(m-2)!}  & \frac{y^{m-1}}{(m-1)!}  &   \frac{y^m}{m!} \\
\frac{y^{m-3}}{(m-3)!}  & \frac{y^{m-2}}{(m-2)!}  &   \frac{y^{m-1}}{(m-1)!}
 \end{array} \right | \frac{(1-x)^k}{k!} dx dy,
$$
then, (5.5) follows from the well-known Beta integral
$$
\int_0^1 x^r (1-x)^s dx = \frac{r! s!}{(r+s+1)!}.
$$

\vspace{5pt}

(3). It suffices to evaluate the following multiple integral
$$
\mathop {\int {\int {...\int {} } } }\limits_{\tiny \begin{array}{*{10}l}
   0 < x_1 < y_1 <  z_1   \\
    \;\;\;\;\;\;\;\; \wedge \;\;\;\;\; \;\;\;\;\; \;\;\; \wedge   \\
    \;\;\;\;\;\;\; x_2  \;\;\;\;\;\;\; < z_2  \\
    \;\;\;\;\;\;\;\; \wedge \;\;\;\; \wedge\;\;\;\; \;\; \wedge   \\
   \;\;\;\;\;\;\; x_3  <  y_2  <  z_3 < 1
\end{array} } \frac{x_3^k}{k!} \; dx_1 \cdots dz_3.
$$
The domain of integration may be decomposed into the following four parts:
\begin{eqnarray*}
D_1 & = & \left \{ \mathop {\tiny \begin{array}{*{10}l}
   0 < x_1 < x_2 <  x_3 < y_1 < y_2 < z_3 < 1   \\
    \;\;\;\;\;\;\;\;\;\;\;\;\;\;\;\;\;\;\;\;\;\;\;\;\;\;\;\;\;\;\;\; \wedge \;\;\;\;\; \;\;\;\;\;\;\; \;\;\; \vee   \\
    \;\;\;\;\;\;\;\;\;\;\;\;\;\;\;\;\;\;\;\;\;\;\;\;\;\;\;\;\;\;\; z_1  \;\;\;\; < \;\;\;\;\; z_2
   \end{array} } \right \},  \;\;\;\;
 D_2  =  \left \{ \mathop {\tiny \begin{array}{*{10}l}
   0 < x_1 <\; y_1   \\
    \;\;\;\;\;\;\;\; \wedge \;\;\;\;\;\;  \wedge   \\
    \;\;\;\;\;\;\; x_2  <  x_3  <   z_1  <  z_2  \\
  \; \; \;\;\;\;\;\;\;\; \; \;\;\;\;\; \wedge \;\;\;\;\;\;\;\;\;\;\; \;\;\; \wedge   \\
    \;\; \; \;\;\;\;\;\;\; \;\;\;\;\;\; y_2 \;\;\;\;< \;\;\;  z_3  < 1
    \end{array}} \right \}, \\
 D_3 & = & \left \{\mathop {\tiny \begin{array}{*{10}l}
\; \; \;\;\;\;\;\;\;\; \; \;\;\;\;\;\; \; \;\;\;\;\; z_2  <  z_3  <  1   \\
    \;\;\;\;\;\;\;\; \; \;\;\;\;\;\;\; \;\;\;\;\;\;\;\; \vee \;\;\;\;\;\;  \vee   \\
0 <  x_1  <  x_2  <   x_3  <  y_2  \\
 \; \; \;\;\;\;\; \wedge \;\;\;\;\;\;\;\;\;\;\; \;\;\;\; \vee   \\
     \; \;\;\;\;\;\; y_1 \;\;\;\;< \;\;\;\;  z_1
    \end{array}} \right \}, \;\;\;\;
D_4  =  \left \{ \mathop {\tiny \begin{array}{*{10}l}
   0 < x_1 < x_2 <  z_2 < x_3 < y_2 < z_3 < 1   \\
   \;\;\;\;\;\;\; \wedge \;\;\;\;\; \;\;\;\;\;\;\; \;\;\; \vee   \\
    \;\;\;\;\;\;\; y_1  \;\;\;\; < \;\;\;\;\; z_1
   \end{array} } \right \},
\end{eqnarray*}
we have
\begin{eqnarray*}
\mathop {\int{\int}} \limits_{\tiny D_1}  & = & \frac{3}{2} \frac{(k+2)!}{k! (k+8)!}, \;\;\;\;\;\;\;  \mathop {\int{\int}} \limits_{\tiny D_2}   =  \frac{(k+3)!}{k! (k+8)!}, \\
\mathop {\int{\int}} \limits_{\tiny D_3}  & = & \frac{1}{4} \frac{(k+4)!}{k! (k+8)!}, \;\;\;\;\;\;\;  \mathop {\int{\int}} \limits_{\tiny D_4}   =  \frac{1}{40} \frac{(k+5)!}{ k! (k+8)!}.
\end{eqnarray*}
Therefore, the proof of theorem 5.2 is complete. \,\,\,\,\, $\Box$

Finally, we note that the order statistics model could  be also used to count the number of solid standard Young tableaux.

\end{document}